\documentclass[11pt]{article}

\usepackage{amsmath,amsthm,amsfonts,amssymb,amsopn,amscd, mathrsfs}
\usepackage{latexsym}
\usepackage{graphics}

\newtheorem{theorem}{Theorem}[section]
\newtheorem{lemma}[theorem]{Lemma}

\newtheorem{corollary}[theorem]{Corollary}

\newtheorem{proposition}[theorem]{Proposition}

\theoremstyle{definition}
\newtheorem{definition}[theorem]{Definition}

\theoremstyle{remark}

\newcommand{\norm}[1]{{\left\| #1 \right\|}}
\newcommand{\defset}[2]{{\left\{ #1 \left| #2 \right. \right\}}}
\newcommand{\corner}[2]{{\mathrm{c}_{ #1 }\left( #2 \right)}}
\newcommand{\offdiag}[2]{{\mathrm{od}_{ #1 }\left( #2 \right)}}

\title{Inclusions and positive cones of von Neumann algebras}
\author{{\bf Yoh Tanimoto} \\
        Graduate School of Mathematical Sciences \\
        University of Tokyo, Komaba, Tokyo, 153-8914, JAPAN \\
        e-mail: {\tt hoyt@ms.u-tokyo.ac.jp} }
\date{}

\begin{document}
\maketitle
\begin{abstract}
We consider cones in a Hilbert space associated to
two von Neumann algebras and determine when one
algebra is included in the other.
If a cone is associated to a von Neumann algebra,
the Jordan structure is naturally recovered
from it and we can characterize projections
of the given von Neumann algebra with the structure
in some special situations.
\end{abstract}

\section{Introduction}
The natural positive cone
$\mathcal{P}^\natural = \overline{\Delta^\frac{1}{4}\mathcal{M}_+\xi_0}$
plays a significant role in the theory of von Neumann algebras
 (see, for example, \cite{araki, connes1})
where $\mathcal{M}$ is a von Neumann algebra,
$\xi_0$ is a cyclic separating vector for $\mathcal{M}$
and $\Delta$ is the Tomita-Takesaki modular operator associated to $\xi_0$.
Among them, the result of Connes \cite{connes2} is
of particular interest which characterized the
natural positive cones with their geometric properties called
selfpolarity, facial homogeneity and orientability, and showed
that if two von Neumann algebras $\mathcal{M}$ and
$\mathcal{N}$ share a same cone, then
there is a central projection $q$ of $\mathcal{M}$ such that
$N = q\mathcal{M} \oplus q^\perp \mathcal{M}^\prime$.
Connes used the Lie algebra with an involution
of the linear transformation group of $\mathcal{P}^\natural$
in his paper.

In the present paper, instead of $\mathcal{P}^\natural$,
we study $\mathcal{P}^\sharp = \overline{\mathcal{M}_+\xi_0}$,
which holds more informations of $\mathcal{M}$, for example,
the subalgebra structure.

In the second section, we study what occurs when
$\overline{\mathcal{N}_+\xi_0} \subset \mathcal{P}^\sharp$ where
$\mathcal{N}$ is another von Neumann algebra.
We consider first the case when $\xi_0$ is not cyclic for $\mathcal{N}$
and then assume the cyclicity.
It turns out that in the latter case $\mathcal{N}$ is included in
$\mathcal{M}$ except the part where $\xi_0$ is tracial.

In the third section, we characterize central projections
of $\mathcal{M}$ in terms of $\mathcal{P}^\sharp$.
A projection $p$ is in $\mathcal{M} \cap \mathcal{M}^\prime$
if and only if $p$ and its orthogonal complement $p^\perp$ preserve
$\mathcal{P}^\sharp$.

In the fourth and fifth sections, the Jordan structure on
$\mathcal{P}^\sharp$ is studied. We can recover the lattice structure of
projections and the operator norm from the order structure of
$\mathcal{P}^\sharp$. Then we can define the square
operation on $\mathcal{P}^\sharp$.

In the final section, using the Jordan structure, a characterization
of projections in $\mathcal{M}$ is obtained when the modular automorphism
with respect to $\xi_0$ acts ergodically.

The result of the second section has an easy application
to the theory of half-sided modular inclusions \cite{wiesbrock, az}.
Let $\{U(t)\}$ be a one-parameter group of unitary operators
with a generator $H$ which kills $\xi_0$.
Assume that $\mathcal{M}$ is a factor of type $\mathrm{III}_1$
(or more generally a properly infinite algebra).
It is easy to see that
$U(t)\mathcal{M}U(t)^* \subset \mathcal{M}$ for $t \ge 0$ if
and only if $U(t)$ preserves $\mathcal{P}^\sharp$ for $t \ge 0$.
A similar result for $\mathcal{P}^\natural$ and $\{e^{-tH}\}$ has been
obtained by Borchers with additional conditions on $H$ \cite{borchers}.

Davidson has obtained conditions for $\{U(t)\}$ to generate a one-parameter
semigroup of endomorphisms \cite{davidson}. The relations with the modular group
have been shown to be important in his study.

\section{Inclusions of positive cones}
Let $\mathcal{M}$ be a von Neumann algebra
acting on a Hilbert space $\mathcal{H}$ and
$\xi_0$ be a cyclic separating vector for $\mathcal{M}$.
We denote the modular group by $\Delta^{it}$, the modular conjugation by $J$,
modular automorphism by $\sigma_t$ and the canonical involution by
$S = J\Delta^\frac{1}{2}$. The positive cone associated to
$\xi_0$ is denoted by $\mathcal{P}^\sharp = \overline{\mathcal{M}_+\xi_0}$.

Suppose there is another von Neumann algebra $\mathcal{N}$ such that
$\overline{\mathcal{N}_+\xi_0} \subset \mathcal{P}^\sharp$.
We can define a positive contractive map $\alpha$ from $\mathcal{N}$
into $\mathcal{M}$ as follows.

\begin{lemma}\label{mapping}
For $a \in \mathcal{N}_+$ there is the unique positive element
$\alpha(a) \in \mathcal{M}$ satisfying $a\xi_0 = \alpha(a)\xi_0$.
In addition, $\alpha$ is contractive on $\mathcal{M}_+$.
\begin{proof}
By the assumption, we have $a\xi_0 \in \mathcal{P}^\sharp$.
Recall that for a vector $a\xi_0$ in $\mathcal{P}^\sharp$ there is a
positive linear operator $\alpha(a)$ affiliated to $\mathcal{M}$ such
that $a\xi_0 = \alpha(a)\xi_0$ \cite{sz}.

Since $\norm{a}I-a$ is positive, we have$(\norm{a}I-a)\xi_0 \in
\mathcal{P}^\sharp$.
This implies, for every $y \in \mathcal{M}^\prime$,
\begin{eqnarray*}
 \langle \alpha(a)y\xi_0, y\xi_0 \rangle
  &=& \langle \alpha(a)\xi_0, y^*y\xi_0 \rangle \\
  &=& \langle a\xi_0, y^*y\xi_0 \rangle \\
  &\le& \norm{a} \langle \xi_0, y^*y\xi_0 \rangle = \norm{a} \norm{y\xi_0}^2.
\end{eqnarray*}
Hence $\alpha(a)$ is bounded and in $\mathcal{M}$. 
\end{proof}
\end{lemma}

We can easily see that $\alpha$ extends to $\mathcal{N}$ by linearity.
Since $\alpha$ is contractive on $\mathcal{N}_+$, $\alpha$ is bounded
on $\mathcal{N}_{sa}$.

\begin{lemma}\label{proj}
The map $\alpha$ maps every projection to a projection.
\begin{proof}
Take a projection $e \in \mathcal{N}$.
Note that, since $\alpha$ maps $\mathcal{N}_+$ into $\mathcal{M}_+$ and
is contractive, we have $\alpha(e) \ge \alpha(e)^2$.

Recall that, by the definition of $\alpha$, we have
$\alpha(e)\xi_0 = e\xi_0$.
We calculate as follows.
\begin{eqnarray*}
 \langle \alpha(e)^2\xi_0,\xi_0 \rangle
 &=& \langle \alpha(e)\xi_0, \alpha(e)\xi_0 \rangle \\
 &=& \langle e\xi_0, e\xi_0 \rangle \\
 &=& \langle e\xi_0, \xi_0 \rangle \\
 &=& \langle \alpha(e)\xi_0, \xi_0 \rangle.
\end{eqnarray*}

This implies that
$\langle \left(\alpha(e)-\alpha(e)^2\right)\xi_0, \xi_0 \rangle = 0$.
As we noted above, $\alpha(e)-\alpha(e)^2$ must be positive,
hence the vector $\left(\alpha(e)-\alpha(e)^2\right)^\frac{1}{2}\xi_0$
must vanish.
By the separating property of $\xi_0$, we see $\alpha(e) = \alpha(e)^2$.
\end{proof}
\end{lemma}

Recall that a linear mapping $\phi$ which preserves every anticommutator is
called a Jordan homomorphism:
\begin{equation*}
 \phi(xy + yx) = \phi(x)\phi(y)+\phi(y)\phi(x).
\end{equation*}
Now we show the following lemma. The proof of it is essentially
taken from \cite{kadison}.

\begin{lemma}
The map $\alpha$ is a Jordan homomorphism.
\begin{proof}
Let $e$ and $f$ be mutually orthogonal projections in $\mathcal{N}$. Then $e+f$,
$\alpha(e)$, $\alpha(f)$ and $\alpha(e)+\alpha(f)$ are projections. We see the
range of $\alpha(e)$ and the range of $\alpha(f)$ are mutually orthogonal because if not, then the sum
$\alpha(e)+\alpha(f)$ could not be a projection. This implies that
\begin{equation*}
\alpha(e)\alpha(f) = \alpha(f)\alpha(e) = 0.
\end{equation*}

In particular, $\alpha$ maps the positive (resp. negative) part of
a self-adjoint element $x$ to the positive (reps. negative) part of
$\alpha(x)$. From this we see that $\alpha$ is contractive
on $\mathcal{N}_{sa}$.

Next suppose we have commuting projections $e,f \in \mathcal{N}$. 
Remark that, since $ef \le e$, positivity of $\alpha$ assures
$\alpha(ef) \le \alpha(e)$. Recalling that in this case $ef$ and $e$ are
projections, we see the range of $\alpha(ef)$ is included in the range
of $\alpha(e)$. Thus we have $\alpha(ef)\alpha(e) = \alpha(ef)$.

Now noting $e-ef$ and $f$ are mutually orthogonal projections, 
we have
\begin{equation*}
  0 = \alpha(e-ef)\alpha(e) = \alpha(e)\alpha(f)-\alpha(ef).
\end{equation*}
Hence $\alpha$ preserves products of commuting projections.

Since every self-adjoint element in a von Neumann algebra is
a uniform limit of linear combinations of mutually orthogonal projections,
and since $\alpha$ is continuous in norm on $\mathcal{N}_{sa}$,
$\alpha$ preserves products of commuting self-adjoint elements.
In particular, $\alpha$ preserves the square of self-adjoint
elements.

This implies that, firstly, $\alpha$ preserves Jordan products
of self-adjoint elements $ab+ba = (a+b)^2-a^2-b^2$.
This shows
\begin{eqnarray*}
 \alpha(ab+ba) &=& \alpha\left((a+b)^2\right)-\alpha(a^2)-\alpha(b^2) \\
 &=& \alpha(a+b)^2 - \alpha(a)^2 - \alpha(b)^2 \\
 &=& \alpha(a)\alpha(b)+\alpha(b)\alpha(a).
\end{eqnarray*}

Secondly, $\alpha$ preserves squares of arbitrary elements
$(a+ib)^2 = a^2 + i(ab+ba) -b^2$:
\begin{eqnarray*}
 \alpha\left((a+ib)^2)\right)
 &=& \alpha\left(a^2 + i(ab+ba) -b^2\right) \\
 &=& \alpha(a^2) + i\alpha(ab+ba) - \alpha(b^2) \\
 &=& \alpha(a)^2 + i\left(\alpha(a)\alpha(b)+\alpha(b)\alpha(a)\right)
      -\alpha(b)^2 \\
 &=& \left(\alpha(a)+i\alpha(b)\right)^2.
\end{eqnarray*}

 Finally, $\alpha$ preserves Jordan products of arbitrary elements
$xy+yx = (x+y)^2-x^2-y^2$:
\begin{eqnarray*}
 \alpha(xy+yx) &=& \alpha\left((x+y)^2\right)-\alpha(x^2)-\alpha(y^2) \\
 &=& \alpha(x+y)^2 - \alpha(x)^2 - \alpha(y)^2 \\
 &=& \alpha(x)\alpha(y)+\alpha(y)\alpha(x).
\end{eqnarray*}
This completes the proof. 
\end{proof}
\end{lemma}

Here we need the following result on Jordan homomorphisms of Jacobson and Rickart \cite{jr}.

\begin{proposition}\label{JR}
Suppose $\phi$ is a unital Jordan homomorphism from an algebra
$\mathcal{A}$ into $\mathcal{B}$. Suppose further that $\mathcal{A}$ has
a system of matrix units. Then there is a central idempotent $g$ of
the algebra generated by $\phi(\mathcal{A})$ such that
$\phi(\cdot)g$ is homomorphic and $\phi(\cdot)(I-g)$
is antihomomorphic.
\end{proposition}

Note that every von Neumann algebra $\mathcal{N}$ decomposes into
the commutative part, the $\mathrm{I}_n$ parts, the $\mathrm{II}_1$ part, and
the properly infinite part. On the first one $\alpha$ causes no problem
and on the remaining parts we can apply Proposition \ref{JR} to the
case in which $\phi = \alpha$,
$\mathcal{A} = \mathcal{N}$, $\mathcal{B} = \mathcal{M}$. Examining the proof,
we see if $\phi$ is self-adjoint, then $g$ is a central projection
of $\alpha(\mathcal{N})^{\prime\prime}$ (the argument here is due to
Kadison \cite{kadison}).

Next, we show the normality of $\alpha$.
\begin{lemma}
The map $\alpha$ is a normal linear mapping from $\mathcal{N}$ into
$\mathcal{M}$.
\begin{proof}
We only have to show that for any normal functional $\varphi$ on $\mathcal{M}$
the functional $\varphi\circ\alpha$ on $\mathcal{N}$ is normal.
Note that, since $\mathcal{M}$ has a separating vector $\xi_0$, we may assume
$\varphi(\cdot) = \langle \cdot\eta_1, \eta_2 \rangle$ for some 
$\eta_1, \eta_2 \in \mathcal{H}$. 

Recall that a linear functional on a von Neumann algebra is normal
if and only if it is continuous on every bounded set in the weak operator
topology.

Now suppose that we have a convergent bounded net in the weak operator topology
$x_i \to x$ in $\mathcal{N}$. Obviously $\{x_i\xi_0\}$ converges to $x\xi_0$
weakly.
By the definition of $\alpha$,
we see $\{\alpha(x_i)\xi_0\}$ converges to $\alpha(x)\xi_0$ weakly.
We have, for any $y_1, y_2 \in \mathcal{M}^\prime$,
\begin{eqnarray*}
  \langle \alpha(x_i)y_1\xi_0, y_2\xi_0 \rangle
   &=& \langle y_1\alpha(x_i)\xi_0, y_2\xi_0 \rangle \\
   &=& \langle \alpha(x_i)\xi_0, y_1^*y_2\xi_0 \rangle \\
   &\to& \langle \alpha(x)\xi_0, y_1^*y_2\xi_0 \rangle \\
   &=& \langle \alpha(x)y_1\xi_0, y_2\xi_0 \rangle.
\end{eqnarray*}

First we assume $\{x_i\}$ is a net of self-adjoint elements.
Then for arbitrary $\eta_1, \eta_2 \in \mathcal{H}$ the convergence
$\langle \alpha(x_i)\eta_1, \eta_2 \rangle
  \to \langle \alpha(x)\eta_1, \eta_2 \rangle$
holds since $\{x_i\}$ is a bounded net, 
$\alpha$ is contractive on $\mathcal{N}_{sa}$, and
$\xi_0$ is cyclic for $\mathcal{M}^\prime$. 

Then we can obtain the convergence for arbitrary bounded
convergent net in WOT $\{x_i\}$ since we have the decomposition
\begin{equation*}
  x_i = \frac{x_i+x_i^*}{2} + i\frac{x_i-x_i^*}{2i}
\end{equation*}
and each part of the net is self-adjoint or antiself-adjoint,
bounded and WOT-converging. 
\end{proof}
\end{lemma}

We combine this lemma and the proposition of Jacobson and Rickart
to get the following.

\begin{lemma}
There is a normal homomorphism $\beta$ and normal antihomomorphism $\gamma$
of $\mathcal{N}$ into $\mathcal{M}$ such that $\alpha(x) = \beta(x) + \gamma(x)$
and the the range of $\beta$ and $\gamma$ are mutually orthogonal.

In addition, there are central projections $e, f \in \mathcal{N}$ and a
central projection $g \in \alpha(\mathcal{N})^{\prime\prime}$ such that
$\alpha(e\ \cdot)g = \beta(\cdot)$ is an isomorphism of $\mathcal{N}e$ and
$\alpha(f\ \cdot)g^\perp = \gamma(\cdot)$ is an antiisomorphism of $\mathcal{N}f$.
\begin{proof}
We know from Proposition \ref{JR} that there is a central projection $g \in
\alpha(\mathcal{N})^{\prime\prime}$ such that $\beta(\cdot) = \alpha(\cdot)g$
is a homomorphism of $\mathcal{N}$ and $\gamma(\cdot) = \alpha(\cdot)g^\perp$ is 
an antihomomorphism of $\mathcal{N}f$.
 Then just take $e$ as the support of $\beta$ and $f$ as the support of
$\gamma$. Since $\alpha$ is normal, so are $\beta$ and $\gamma$ and
the definitions of $e$ and $f$ are legitimate. 
\end{proof}
\end{lemma}

\begin{lemma}\label{prop_inf}
The von Neumann algebra $\mathcal{N}f$ is finite.
\begin{proof}
Let $\mathcal{N}h$ be the properly infinite part of $\mathcal{N}f$.
We have
$g^\perp\alpha(xy) = g^\perp\alpha(y)\alpha(x) = \alpha(y)g^\perp\alpha(x)$
for $x,y \in \mathcal{N}h$.

Again take $x,y \in \mathcal{N}h$. By the definition of $\alpha$,
we have
\begin{eqnarray*}
  g^\perp xy\xi_0 &=& g^\perp \alpha(xy)\xi_0 \\
               &=& \alpha(y)g^\perp \alpha(x)\xi_0 \\
  \left\langle g^\perp xy\xi_0, \xi_0 \right\rangle 
   &=& \left\langle \alpha(y) g^\perp \alpha(x)\xi_0, \xi_0 \right\rangle \\
   &=& \left\langle  g^\perp \alpha(x)\xi_0, \alpha(y^*)\xi_0 \right\rangle \\
   &=& \left\langle  g^\perp x\xi_0, y^*\xi_0 \right\rangle \\
   &=& \left\langle y g^\perp x\xi_0, \xi_0 \right\rangle.
\end{eqnarray*}

Since $\mathcal{N}h$ is properly infinite, there is a sequence of isometries
$\{v_n\} \subset \mathcal{N}h$ such that $v_nv_n^* \to 0$ in SOT-topology
(That they are isometries means $v_n^*v_n = h$). Now
\begin{eqnarray*}
 \left\langle \gamma(h)\xi_0, \xi_0 \right\rangle
  &=& \left\langle  g^\perp h\xi_0, \xi_0 \right\rangle \\
  &=& \left\langle  g^\perp v_n^*v_n\xi_0, \xi_0 \right\rangle \\
  &=& \left\langle v_n g^\perp v_n^*\xi_0, \xi_0 \right\rangle \\
  &\le& \left\langle v_nv_n^*\xi_0, \xi_0 \right\rangle \to 0.
\end{eqnarray*}
But since $\gamma(h)$ is a projection in
$\alpha(\mathcal{N})^{\prime\prime} \subset \mathcal{M}$
and since $\xi_0$ is separating for $\mathcal{M}$, $\gamma(h)$ must be zero.
Recalling that $h$ is a subprojection of $f$ and that $f$ is the support of
$\gamma$, we see that $h = 0$.
\end{proof}
\end{lemma}

\begin{theorem}\label{gen}
Let $\mathcal{M}$ and $\mathcal{N}$ be von Neumann algebras and $\xi_0$ is a
cyclic separating vector for $\mathcal{M}$.
Suppose $\overline{\mathcal{N}_+\xi_0} \subset \mathcal{P}^\sharp$.

Then we have two disjoint possibilities:
\begin{enumerate}
\item The von Neumann algebra $\mathcal{M}$ has a subalgebra $\mathcal{M}_1$
      such that $\overline{\mathcal{M}_{1+}\xi_0} =
      \overline{\mathcal{N}_+\xi_0}$.
\item For any subalgebra $\mathcal{M}_2$ of $\mathcal{M}$,
      its ``sharpened cone''
      $\overline{\mathcal{M}_{2+}\xi_0}$ cannot coincide with
      $\overline{\mathcal{N}_+\xi_0}$ and $\mathcal{N}$ has a finite ideal
      $\mathcal{N}_1$ such that
      there is a subalgebra of $\mathcal{M}$ which is isomorphic
      to the direct sum of $\mathcal{N}_1$ and $\mathcal{N}_1^{\mathrm{opp}}$.
\end{enumerate}
\begin{proof}
Suppose that $e$ and $f$ defined above are mutually orthogonal.
Then let us define $\mathcal{M}_1 = \alpha(\mathcal{N})$.
Since we have $ef = 0$, it decomposes as follows.
\begin{eqnarray*}
 \alpha(\mathcal{N}) &=& \alpha\left(\mathcal{N}[e+ e^\perp][f+ f^\perp ]\right) \\
 &=& \alpha\left(\mathcal{N}[e f^\perp +f e^\perp+ e^\perp f^\perp ]\right) \\
 &=& \beta\left(\mathcal{N}e f^\perp \right)+\gamma\left(\mathcal{N}f e^\perp \right),
\end{eqnarray*}
by noting that $\mathcal{N} e^\perp  f^\perp $ is the kernel of $\alpha$.

Since the range of $\beta$ and $\gamma$ are mutually orthogonal, and since
$e$ and $f$ are central projections,
$\alpha(\mathcal{N})$ is a direct sum of $\beta\left(\mathcal{N}e f^\perp \right)$
and $\gamma\left(\mathcal{N}f e^\perp \right)$.

Let $a$ be a positive element of $\mathcal{N}$. Then we have
\begin{eqnarray*}
 a\xi_0 &=& \alpha(a)\xi_0 \\
        &=& \beta(ae)\xi_0 + \gamma(af)\xi_0 \\
        &=& \beta(ae f^\perp )\xi_0 + \gamma(af e^\perp )\xi_0.
\end{eqnarray*}
Conversely it is easy to see that for $b \in \alpha(\mathcal{N})_+$
there is $a \in \mathcal{N}_+$ such that $\alpha(a) = b$, hence we have
$a\xi_0 = b\xi_0$. This completes the proof of the claimed equality
$\overline{\mathcal{M}_{1+}\xi_0} = \overline{\mathcal{N}+\xi_0}$.

Next, we assume that $ef \neq 0$. Note that $\mathcal{N}ef$ is
noncommutative since by the definition of $\beta$ and $\gamma$
the commutative part of $\mathcal{N}$ is left to $\beta$.
In particular $g$ is a nontrivial central projection
in $\alpha(\mathcal{N}ef)^{\prime\prime}$.
By Lemma \ref{prop_inf}, $\mathcal{N}ef$ is finite.
One can easily see that $\alpha(\mathcal{N}ef)^{\prime\prime}$ is a subalgebra
of $\mathcal{M}$ which decomposes into the direct sum of
$\beta(\mathcal{N}ef)$ and $\gamma(\mathcal{N}ef)$ where
the latter is isomorphic to $(\mathcal{N}ef)^{\mathrm{opp}}$.

What remains to prove is that for any subalgebra $\mathcal{M}_2$ of
$\mathcal{M}$ we cannot have the equality
$\overline{(\mathcal{N}ef)_+\xi_0} = \overline{\mathcal{M}_{2+}\xi_0}$.
To see this impossibility, recall that
\begin{equation*}
 \overline{\mathcal{M}_+\xi_0} = \defset{A\xi_0}{ A \mbox{ is a closed
 positive operator affiliated to }\mathcal{M}},
\end{equation*}
since $\xi_0$ is a separating vector for $\mathcal{M}$ \cite{sz}.
Similarly we have
\begin{equation*}
 \overline{\mathcal{M}_{2+}\xi_0} = \defset{A\xi_0}{A \mbox{ is a closed
 positive operator affiliated to }\mathcal{M}_2}.
\end{equation*}

Now suppose $a\xi_0 \in \overline{\mathcal{M}_{2+}\xi_0}$ for a positive element
$a$ of $\mathcal{N}ef$. By the above remark, we have a positive operator
$A$ affiliated to $\mathcal{M}_2$ such that $a\xi_0 = \alpha(a)\xi_0 = A\xi_0$. 
Then for $y \in \mathcal{M}^\prime$ we have
\begin{equation*}
  \alpha(a)y\xi_0 = y\alpha(a)\xi_0 = yA\xi_0 = Ay\xi_0,
\end{equation*}
hence $A$ is bounded and $\alpha(a) = A$.
This implies $\alpha(a) \in \mathcal{M}_2$ and
$\alpha(\mathcal{N}ef) \subset \mathcal{M}_2$. But by Proposition \ref{JR}
$\alpha(\mathcal{N}ef)$ generates
$\beta(\mathcal{N}ef) \oplus \gamma(\mathcal{N}ef)$. We have
$\beta(\mathcal{N}ef) \oplus \gamma(\mathcal{N}ef) \subset \mathcal{M}_2$.

We will show that this leads to a contradiction. 
By the observation above we see that $\overline{\mathcal{M}_{2+}\xi_0}$ contains
vectors of the form $ga\xi_0,  g^\perp b\xi_0$ where $a, b \in (\mathcal{N}ef)_+$.

Suppose the contrary that $ga\xi_0 \in \overline{(\mathcal{N}ef)_+\xi_0}$.
By the argument similar to the above one,
there is a self-adjoint positive operator $A$
affiliated to $\mathcal{N}ef$ such that $A\xi_0 = ga\xi_0$.
Then $ g^\perp A\xi_0 = 0$. Noting that $f$ is the support of $\gamma$ and
that $\xi_0$ is separating for $\mathcal{M}$, we see
$ g^\perp e_A\xi_0 = \gamma(e_A)\xi_0$ cannot vanish for any nontrivial projection
$e_A$ of $\mathcal{N}ef$.

There are a spectral projection $e_A$ of $A$, a positive scalar $\epsilon$
and $y \in \mathcal{M}^\prime$ such that $A \ge \epsilon e_A$ and
$\langle \gamma(e_A)y\xi_0, y\xi_0 \rangle > 0$. 
Remark that 
\begin{eqnarray*}
   g^\perp (A-\epsilon e_A)\xi_0 &\in&  g^\perp \overline{(\mathcal{N}ef)_+\xi_0} \\
    &\subset& \overline{ g^\perp (\mathcal{N}ef)_+\xi_0} \\
    &=& \overline{\gamma(\mathcal{N}ef)_+\xi_0}.
\end{eqnarray*}

Then we have 
\begin{eqnarray*}
  0 &=& \langle y g^\perp A\xi_0, y\xi_0 \rangle \\
    &=& \langle  g^\perp A\xi_0, y^*y\xi_0 \rangle \\
    &=& \langle  g^\perp (A-\epsilon e_A)\xi_0, y^*y\xi_0 \rangle 
        + \langle  g^\perp \epsilon e_A\xi_0, y^*y\xi_0 \rangle \\
    &\ge& \langle  g^\perp \epsilon e_A\xi_0, y^*y\xi_0 \rangle \\
    &=& \langle y\gamma(\epsilon e_A)\xi_0, y\xi_0 \rangle \\
    &=& \epsilon \langle \gamma(e_A)y\xi_0, y\xi_0 \rangle \\
    &>& 0.
\end{eqnarray*}

This contradiction completes the proof of that
$\overline{(\mathcal{N}ef)_+\xi_0} \neq \overline{\mathcal{M}_{2+}\xi_0}$.
\end{proof}
\end{theorem}

If we further assume the cyclicity of $\xi_0$ for $\mathcal{N}$, we have a
stronger result. For the proof of it, we need the following lemma.
This can be found, for example in \cite{baumgaertel}, but here we present
another simple proof.
\begin{lemma}\label{fin_coinc}
If $\mathcal{A} \subset \mathcal{B}$ is a proper inclusion of von Neumann algebras
on a Hilbert space $\mathcal{K}$ 
and if $\zeta$ is a common cyclic separating vector,
then $\mathcal{B}$ cannot be finite.
\begin{proof}
Suppose the contrary, that $\mathcal{B}$ is finite.
Then $\mathcal{A}$ must be finite, too.
Hence there is a faithful trace $\tau$ on $\mathcal{B}$. Since
$\zeta$ is separating for $\mathcal{B}$, 
there is a vector $\eta$ such that $\tau(x) = \langle x\eta, \eta \rangle$
by the Radon-Nikodym type theorem.
Since $\tau$ is faithful, $\eta$ must be separating for $\mathcal{B}$.

We can see that $\eta$ is cyclic for $\mathcal{B}$ as follows.
Denote the orthogonal projection onto $\overline{\mathcal{B}\eta}$ by $p$.
By separation verified above,
we have $\overline{\mathcal{B}^\prime\eta} = \mathcal{K}$.
On the other hand, by assumption,
$\overline{\mathcal{B}\zeta} = \overline{\mathcal{B}^\prime\zeta} = \mathcal{K}$.
By the general theory of equivalence of projections,
$p \sim I$ in $\mathcal{B}$. But recalling that $\mathcal{B}$ is finite,
we see that $p = I$, i.e., $\eta$ is cyclic.

By the same reasoning, $\eta$ is cyclic separating tracial for $\mathcal{A}$.
Then the modular conjugations $J_\mathcal{A}$ and $J_\mathcal{B}$ with
respect to $\eta$ must
coincide and we have the required equation.
\begin{equation*}
  \mathcal{A}^\prime \supset \mathcal{B}^\prime
   = J_\mathcal{B}\mathcal{B}J_\mathcal{B} = J_\mathcal{A}\mathcal{B}J_\mathcal{A}
   \supset J_\mathcal{A}\mathcal{A}J_\mathcal{A} = \mathcal{A}^\prime.
\end{equation*}
This contradicts the assumption that the inclusion $\mathcal{A} \subset \mathcal{B}$
is proper. 
\end{proof}
\end{lemma}

\begin{theorem}
Let $\mathcal{M}$ and $\mathcal{N}$ be von Neumann algebras and
$\xi_0$ be a vector cyclic separating for $\mathcal{M}$ and 
cyclic for $\mathcal{N}$.
Suppose $\overline{\mathcal{N}_+\xi_0} \subset \mathcal{P}^\sharp$.

Then we have the following.
\begin{enumerate}
\item The vector $\xi_0$ is also separating for $\mathcal{N}$.
\item There is a central projection $e$ in $\mathcal{N}$ such that 
      $\mathcal{N}e \subset \mathcal{M}$.
\item The vector $e^\perp \xi_0$ is tracial for $\mathcal{N}e^\perp$. 
\item $J_{e^\perp}\mathcal{N}e^\perp J_{e^\perp} \subset
       \mathcal{M}$.
\end{enumerate}

In particular, $\mathcal{N}$ and $\mathcal{N}e \oplus J_{e^\perp}\mathcal{N}e^\perp J_{e^\perp}$
share the same positive cone $\mathcal{P}^\sharp_\mathcal{N}$ where
$\mathcal{N}e \oplus J_{e^\perp}\mathcal{N}e^\perp J_{e^\perp} \subset \mathcal{M}$.
\begin{proof}
First we show that the induction by $g$ realizes $\beta(\cdot) = g\alpha(\cdot)$.
For arbitrary $x,y \in \mathcal{N}$ we have
\begin{eqnarray*}
 gxy\xi_0 &=& g\alpha(xy)\xi_0 \\
          &=& g\alpha(x)\alpha(y)\xi_0 \\
          &=& g\alpha(x)y\xi_0 \\
          &=& \alpha(x)gy\xi_0.
\end{eqnarray*}
Taking it into consideration that $\xi_0$ is cyclic for $\mathcal{N}$,
we see that $gx = g\alpha(x) = \alpha(x)g$. But, since this holds for
arbitrary $x \in \mathcal{N}$, in particular for self-adjoint elements.
If $x = x^*$, then we have
\begin{equation*}
 gx = \alpha(x)g = \left(g\alpha(x)\right)^* = (gx)^* = xg.
\end{equation*}
Since this equation is linear for $x$, we see that
$g \in \mathcal{N}^\prime$ and $gx = g\alpha(x)$.

Now recall that we have decomposed $\alpha$ into a normal homomorphism
$\beta$ and a normal antihomomorphism $\gamma$. We again denote
the support of $\beta$ by $e$ and the support of $\gamma$ by $f$.

Let $\mathcal{N}h$ be the properly infinite part. By Lemma \ref{prop_inf}
the intersection of $h$ and $f$ is trivial. Thus we have
\begin{equation*}
  ghx\xi_0 = hg\alpha(hx)\xi_0 = h\alpha(hx)\xi_0 = hx\xi_0,
\end{equation*}
for $x \in \mathcal{N}$. Cyclicity of $\xi_0$ tells us that $gh = h$.
Then for $hx \in \mathcal{N}h$ we get that
\begin{equation*}
  \alpha(hx) = ghx = hx.
\end{equation*}
In other words, $\alpha$ maps identically on $\mathcal{N}h$.
In particular, $\alpha$ is decomposed by $h$, that is, we have
\begin{equation*}
  h\alpha(h^\perp) = \alpha(h)\alpha(h^\perp) = 0,
\end{equation*}
since $\alpha$ maps orthogonal projections to orthogonal projections.

Note that $h\xi_0$ is cyclic for $\mathcal{N}h$ since $\xi_0$ is cyclic
for $\mathcal{N}$. The vector $h\xi_0$ is also separating for $\mathcal{N}h$
since
\begin{equation*}
  \mathcal{N}h = \alpha(\mathcal{N}h) \subset \mathcal{M}
\end{equation*}
and $\xi_0$ is separating for $\mathcal{M}$.

For the proof of remaining part of the theorem, we may assume $\mathcal{N}$
is finite.

Recall that $ g^\perp $ commutes with $\mathcal{N}$.
Take $x,y \in \mathcal{N}$ and let us calculate
\begin{eqnarray*}
  \langle xyg^\perp \xi_0, g^\perp \xi_0 \rangle
   &=& \langle g^\perp y\xi_0, g^\perp x^*\xi_0 \rangle \\
   &=& \langle g^\perp \alpha(y)\xi_0, g^\perp \alpha(x^*)\xi_0 \rangle \\
   &=& \langle g^\perp \alpha(x)\alpha(y)\xi_0, g^\perp \xi_0 \rangle \\
   &=& \langle g^\perp \alpha(yx)\xi_0, g^\perp \xi_0 \rangle \\
   &=& \langle g^\perp yx\xi_0, g^\perp \xi_0 \rangle \\
   &=& \langle yxg^\perp \xi_0, g^\perp \xi_0 \rangle
\end{eqnarray*}
This shows that $g^\perp \xi_0$ is a tracial vector for $\mathcal{N}g^\perp $.
By assumption, $\xi_0$ is cyclic for $\mathcal{N}$,
hence $g^\perp \xi_0$ is cyclic for $\mathcal{N}g^\perp$.
In addition, it is also separating as follows. If $xg^\perp \xi_0 = 0$ for some
$x \in \mathcal{N}g^\perp $, then
for any $y \in \mathcal{N}g^\perp $ we have
\begin{eqnarray*}
 \norm{xyg^\perp\xi_0}^2 &=& \langle y^*x^*xyg^\perp\xi_0, g^\perp\xi_0 \rangle \\
                &=& \langle xyy^*x^*g^\perp \xi_0, g^\perp \xi_0 \rangle \\
              &\le& \norm{y}^2 \langle xx^*g^\perp \xi_0, g^\perp \xi_0 \rangle \\
                &=& \norm{y}^2 \langle x^*xg^\perp \xi_0, g^\perp \xi_0 \rangle \\
                &=& 0,
\end{eqnarray*}
then the cyclicity implies the separation by $g^\perp \xi_0$.

Now $\mathcal{N}g^\perp $ has the canonical conjugation $J_{g^\perp}$
defined as (the closure of)
\begin{equation*}
 J_{g^\perp}:g^\perp \mathcal{H} \ni x\xi_0 \longmapsto x^*\xi_0
   \in g^\perp \mathcal{H}.
\end{equation*}
On $\mathcal{N}g^\perp $ we have the canonical antihomomorphism
\begin{equation*}
 \mathcal{N}g^\perp \ni x \longmapsto J_{g^\perp}x^*J_{g^\perp}
   \in \mathcal{N}g^\perp .
\end{equation*}

In our situation the composition of the induction by $g^\perp$ and
this antihomomorphism coincide with the composition of $\alpha$ and
the induction by $g^\perp $.
In fact, for any elements $x,y,z \in \mathcal{N}g^\perp$ we have
\begin{eqnarray*}
 \langle J_{g^\perp}(xg^\perp)^*g^\perp J_{g^\perp} yg^\perp \xi_, zg^\perp \xi_0 \rangle 
  &=& \langle z^*g^\perp \xi_0, x^*y^*g^\perp \xi_0 \rangle \\
  &=& \langle yxz^*g^\perp \xi_0, g^\perp \xi_0 \rangle \\
  &=& \langle z^*yxg^\perp \xi_0, g^\perp \xi_0 \rangle \\
  &=& \langle g^\perp yx\xi_0, zg^\perp \xi_0 \rangle \\
  &=& \langle g^\perp \alpha(yx)\xi_0, zg^\perp \xi_0 \rangle \\
  &=& \langle g^\perp \alpha(x)\alpha(y)\xi_0, zg^\perp \xi_0 \rangle \\
  &=& \langle g^\perp \alpha(x)y\xi_0, zg^\perp \xi_0 \rangle \\
  &=& \langle g^\perp \alpha(x)y\xi_0, zg^\perp \xi_0 \rangle.
\end{eqnarray*}
The cyclicity of $g^\perp \xi_0$ shows that $g^\perp \alpha(x) = J_{g^\perp}(xg^\perp)^*J_{g^\perp}$.

Summing up, we get the following formula for $\alpha$:
\begin{eqnarray*}
  \alpha(x) &=& g\alpha(x) +  g^\perp \alpha(x) \\
            &=& gx + J_{g^\perp}g^\perp x^* J_{g^\perp}.
\end{eqnarray*}

Note that $g\xi_0$ is cyclic separating for $\mathcal{N}g$.
In fact, the cyclicity comes from the assumption of $\xi_0$'s 
cyclicity and separating property can be seen by observing
\begin{equation*}
  \mathcal{N}g = g\alpha(\mathcal{N}) \subset \mathcal{M}
\end{equation*}
and by separating property of $\xi_0$ for $\mathcal{M}$.

On the other hand, we have seen that $g^\perp \xi_0$ is cyclic
separating for $\mathcal{N}g^\perp $ in the way proving that
$g^\perp \xi_0$ is a faithful tracial vector.

The direct sum of $\mathcal{N}g$ and
$\mathcal{N}g^\perp $ has a cyclic separating vector $\xi_0$.
These summands are finite because we are assuming
that $\mathcal{N}$ is finite and they are induced 
part of it. Hence $\mathcal{N}g \oplus \mathcal{N}g^\perp $
is also finite.

Clearly $\mathcal{N}$ is a subalgebra of
$\mathcal{N}g \oplus \mathcal{N}g^\perp $. So $\xi_0$ is separating
for $\mathcal{N}$. This is the first statement of the theorem.

Now we have an inclusion of finite von Neumann algebras
\begin{equation*}
  \mathcal{N} \subset \mathcal{N}g \oplus \mathcal{N}g^\perp
\end{equation*}
and $\xi_0$ is a common cyclic separating vector.
Then they must coincide by Lemma \ref{fin_coinc}. This happens only if
$g$ is a projection of $\mathcal{N}$ from the
beginning, i.e, $g$ is a central projection of
$\mathcal{N}$.

Recall that induction by $g$ coincides with the homomorphic
part of $\alpha$. Now we know that $g$ is central.
Then the support $e$ of the homomorphic part $\beta$ must be
exactly $g$.

On the other hand, the intersection $ e^\perp  f^\perp $ of kernels
of the homomorphic part $\beta$ and the antihomomorphic part
$\gamma$ must be trivial. To see this, take $x \in \mathcal{N}$.
We have
\begin{eqnarray*}
   e^\perp  f^\perp x\xi_0 &=& x e^\perp  f^\perp \xi_0 \\
                   &=& x\alpha\left( e^\perp  f^\perp \right)\xi_0 \\
                   &=& 0.
\end{eqnarray*}
Since $\xi_0$ is cyclic for $\mathcal{N}$, we get that
$ e^\perp  f^\perp  = 0$.

Since the induction by $e$ realizes the homomorphic part $\beta$
of $\alpha$, for the antihomomorphic part $\gamma$ it holds
\begin{equation*}
  \gamma(e) =  e^\perp \alpha(e) = \alpha(e)-e\alpha(e) = 0.
\end{equation*}
This implies $e$ must be orthogonal to $f$, which is the
support of $\gamma$. As their intersection vanishes, we get
$f = I-e$.

Recalling $g = e$, we saw that $e^\perp\xi_0$ is a cyclic separating tracial
vector for $\mathcal{N}e^\perp$ and the canonical antiisomorphism with respect
to $e^\perp\xi_0$ coincides with $e^\perp\alpha$. Then the proof of all the statements
in the theorem is done.
\end{proof}
\end{theorem}

\section{Recovery of central projections}
In the following sections we turn to the study of single von Neumann algebra.
Again let $\mathcal{M}$ be a von Neumann algebra and $\xi_0$ be a cyclic separating vector for $\mathcal{M}$.
By Connes' result, $\mathcal{P}^\natural$ determines
$\mathcal{M}$ up to center.

Here we show that the center is easily recovered from $\mathcal{P}^\sharp$.
Let $p$ be a projection $\mathcal{B}(\mathcal{H})$ such that
$p\mathcal{P}^\sharp \subset \mathcal{P}$ and
$ p^\perp \mathcal{P}^\sharp \subset \mathcal{P}^\sharp$.

In this situation, we can define a mapping from $\mathcal{M}$ into
$\mathcal{M}$ using $p$.

\begin{lemma}\label{def_cent}
For every $a \in \mathcal{M}_+$ there is $\alpha(a) \in \mathcal{M}_+$ such that
$pa\xi_0 = \alpha(a)\xi_0$.
\begin{proof}
As in the proof of Lemma \ref{mapping}, we have a positive operator $\alpha(a)$
affiliated to $\mathcal{M}$ such that $pa\xi_0 = \alpha(a)\xi_0$ since $pa\xi_0$
is a vector of the positive cone $\mathcal{P}^\sharp$. This is again bounded
for a different reason. In fact, for $y \in \mathcal{M}^\prime$ we have
\begin{eqnarray*}
\langle \alpha(a)y\xi_0, y\xi_0 \rangle
 &=& \langle \alpha(a)\xi_0, y^*y\xi_0 \rangle \\
 &=& \langle pa\xi_0, y^*y\xi_0 \rangle \\
 &\le& \langle pa\xi_0, y^*y\xi_0 \rangle
   + \langle  p^\perp a\xi_0, y^*y\xi_0 \rangle \\
 &=& \langle a\xi_0, y^*y\xi_0 \rangle \\
 &=& \langle ay\xi_0, y\xi_0 \rangle \\
 &\le& \norm{a} \norm{y\xi_0}^2,
\end{eqnarray*}
where we have used the assumption that $ p^\perp $ preserves $\mathcal{P}^\sharp$.
\end{proof}
\end{lemma}
 From this we see that $\alpha(a) \le a$ as self-adjoint operators.
The map $\alpha$ extends to a linear mapping of $\mathcal{M}$.

\begin{lemma}
The map $\alpha$ maps every projection to a projection.
\begin{proof}
Let $e$ be a projection of $\mathcal{M}$. By the observation above,
we have $\alpha(e) \le e$. Then using the fact $e\alpha(e) = \alpha(e)$
we can calculate
\begin{eqnarray*}
\langle \alpha(e)^2\xi_0 \xi_0 \rangle 
 &=& \langle \alpha(e)\xi_0, \alpha(e)\xi_0 \rangle \\
 &=& \langle pe\xi_0, pe\xi_0 \rangle \\
 &=& \langle pe\xi_0, e\xi_0 \rangle \\
 &=& \langle \alpha(e), e\xi_0 \rangle \\
 &=& \langle \alpha(e), \xi_0 \rangle.
\end{eqnarray*}
We can see that $\alpha(e)^2 = \alpha(e)$ as in the proof of
Lemma \ref{proj}. 
\end{proof}
\end{lemma}

Then the mapping $\alpha$ is a normal Jordan homomorphism and there is a
central projection $g$ of
$\alpha(\mathcal{M})^{\prime\prime}\ \subset \mathcal{M}$ such that
$\alpha(\cdot)g$ is homomorphic and $\alpha(\cdot) g^\perp $ is antihomomorphic.
The proof is similar to the one for the case of subcones.

Now we have the following.
\begin{theorem}
Let $\mathcal{M}$ be a von Neumann algebra acting on a Hilbert
space $\mathcal{H}$, $\xi_0$ be a cyclic separating
vector for $\mathcal{M}$ and
$\mathcal{P}^\sharp = \overline{\mathcal{M}_+\xi_0}$.
Then a projection $p \in \mathcal{B}(\mathcal{H})$ is a
central projection of $\mathcal{M}$ if and only if
$p$ and $p^\perp$ preserve $\mathcal{P}^\sharp$.
\begin{proof}
The ``only if'' part is trivial.

Let $p$ be a projection which and whose orthogonal complement preserve $\mathcal{P}^\sharp$.
Note that $\alpha(x) \in \mathcal{M}$ and that
$\alpha\left(\alpha(x)\right) = \alpha(x)$ holds.
In fact, we have
\begin{equation*}
 \alpha\left(\alpha(x)\right)\xi_0 = p\alpha(x)\xi_0
  = ppx\xi_0 = px\xi_0 = \alpha(x)\xi_0,
\end{equation*}
since $p$ is a projection.

As in the situation of subcones, $\alpha$ is a sum of
a normal homomorphism and a normal antihomomorphism whose
ranges are mutually orthogonal. The kernels of the homomorphism
and the antihomomorphism are central projections of $\mathcal{M}$.
Thus the support of $\alpha$ is the orthogonal complement of the
intersection of these kernels. In particular it is
a central projection $e \in \mathcal{M}$.

Recall that $\alpha(e) \le e$. Take an arbitrary positive element
$a$ from $\mathcal{M}$. If we apply $\alpha$ to $ea-\alpha(ea)$,
since the composition of $\alpha$ and $\alpha$ equals $\alpha$ itself,
we have
\begin{equation*}
 \alpha\left(ea-\alpha(ea)\right) = \alpha(ea) - \alpha(ea) = 0.
\end{equation*}
The argument of the left hand side is less than the support of $\alpha$,
hence it must vanish.
Thus we see that $ea$ is fixed by $\alpha$.
By linearity, this holds for arbitrary element $x \in \mathcal{M}$ instead
of positive element $a$.

Again since $e$ is the support of $\alpha$, we have
$\alpha(x) = \alpha(xe) = xe$.
Comparing this with the definition of $\alpha$ we can
determine $p$.
\begin{eqnarray*}
px\xi_0 &=& \alpha(x)\xi_0 \\
        &=& ex\xi_0
\end{eqnarray*}
With the cyclicity of $\xi_0$ we see that $p$ equals $e$.
In particular, $p$ must be a central projection of $\mathcal{M}$.
\end{proof}
\end{theorem}

\section{Properties of $(\mathcal{P}^\sharp, \xi_0)$}\label{cone}
In this section, we study the properties of $\mathcal{P}^\sharp$ coupled
with a specified vector $\xi_0$. We begin with the following lemma.

Let us write $\zeta \le \eta$ if $\eta - \zeta \in \mathcal{P}^\sharp$.

\begin{lemma}\label{fund}
Let $\zeta$ be a vector in $\mathcal{P}^\sharp$. Then the following hold.
\begin{enumerate}
\item If $\zeta \le \xi_0$, then there is a positive contractive operator
      $a \in \mathcal{M}$ such that $\zeta = a\xi_0$.
      In this case we say that $\zeta$ is contractive.
\item If $\zeta$ is contractive and if $\zeta \perp (\xi_0-\zeta)$,
      then there is a projection $e \in \mathcal{M}$ such that $\zeta = e\xi_0$.
      When these conditions hold, we call $\zeta$ a projective vector.
\item If $\eta$ and $\zeta$ are projective and $\zeta \le \xi_0-\eta$, then
      $e$ and $f$ are mutually orthogonal projections where $\eta = e\xi_0$ and
      $\zeta = f\xi_0$. We say $\eta$ and $\zeta$ are mutually operationally
      orthogonal.
\end{enumerate}
\begin{proof}
The proofs of the first and the second statements are same as in the proofs
of Lemma \ref{mapping} and \ref{proj} respectively. We do not repeat them here.

Suppose $\eta = e\xi_0$, $\zeta = f\xi_0$ and $\eta \le \xi_0-\zeta$.
Then according to this order, $e \le I-f$. When $e$ and $f$ are
projections, this shows the mutual orthogonality.  
\end{proof}
\end{lemma}

We denote the set of contractive vectors by $\mathcal{P}^\sharp_1$.
By the Lemma above, to each vector in $\mathcal{P}^\sharp_1$ there corresponds
a positive contractive operator of $\mathcal{M}$.

Similarly to every vector $\zeta$ in $\mathbb{R}_+\mathcal{P}^\sharp_1$
there corresponds a bounded positive operator $a$ of $\mathcal{M}$.
Put $\mathcal{P}^\sharp_b = \mathbb{R}_+\mathcal{P}^\sharp_1$ and
$\mathcal{K} = \mathbb{R}\mathcal{P}^\sharp_1$.

\begin{lemma}
For an arbitrary vector $\zeta$ in $\mathcal{P}^\sharp_1$ there is
a least projective vector such that $\eta \ge \zeta$.
Let us call $\eta$ the support of $\zeta$.
\begin{proof}
As noted above, there is a positive operator $a$ such that
$\zeta = a\xi_0$.  As we have seen, the order
structure of $\mathcal{P}^\sharp_1$ is consistent with this correspondence.
Let $e$ be the support projection of $a$. Then
we have $\eta = e\xi_0 \ge a\xi_0 = \zeta$.
Hence $\eta$ is the least projective vector in $\mathcal{P}^\sharp_1$. 
\end{proof}
\end{lemma}

\begin{lemma}
Every vector $\zeta$ in $\mathcal{K}$ is uniquely decomposed as
$\zeta = \zeta_+ - \zeta_-$ where $\zeta_+$ and $\zeta_-$ are vectors of
$\mathcal{P}^\sharp_b$ and supports of $\zeta_+$ and $\zeta_-$ are mutually
operationally orthogonal.
\begin{proof}
Since every vector in $\mathcal{P}^\sharp_1$ corresponds to
a positive contractive operator in $\mathcal{M}$, vectors of
$\mathcal{P}^\sharp_b$ (resp. $\mathcal{K}$) correspond to
positive operators (resp. self-adjoint operators).

Now the lemma follows from the theory of self-adjoint operators.
The self-adjoint operator $z$ corresponding to $\zeta$ has the
Jordan decomposition $z = z_+ - z_-$ where $z_+$ and $z_-$ are positive
operators of $\mathcal{M}$ whose supports are mutually orthogonal.
By Lemma \ref{fund}, $\zeta$ has the corresponding decomposition. 
\end{proof}
\end{lemma}

\begin{lemma}
The cone $\mathcal{P}^\sharp_b$ is dense in $\mathcal{P}^\sharp$.
\begin{proof}
For each vector $\zeta$ in $\mathcal{P}^\sharp$ there is a positive self-adjoint
linear operator $A$ affiliated to $\mathcal{M}$ such that
$\zeta = A\xi_0$\cite{sz}.
Let $E_A$ be the spectral measure associated to $A$.
Then $AE_A\left([0,n]\right)$ is bounded positive operator in $\mathcal{M}$.
It is well known that $\left\{AE_A\left([0,n]\right)\xi_0\right\}$ converges to
$A\xi_0$. 
\end{proof}
\end{lemma}

In addition, we can recover the operator norm in terms of
$\mathcal{P}^\sharp_b$.
For $\zeta \in \mathcal{P}^\sharp_b$ we define the new ``sharp'' norm
$\norm{\zeta}_\sharp$ as follows.
\begin{equation*}
  \norm{\zeta}_\sharp
   = \mathrm{sup}\defset{c \ge 0}{ \frac{1}{c}\zeta \le \xi_0}.
\end{equation*}

\begin{lemma}
If $a \in \mathcal{M}_+$ and $\zeta = a\xi_0$, then $\norm{\zeta}_\sharp =
\norm{a}$.
\begin{proof}
We only have to note that $ca\xi_0 \le \xi_0$ if and only if $ca \le I$.
Then the spectral decomposition of $a$ completes the proof. 
\end{proof}
\end{lemma}

The set $\mathcal{K}$ is a real linear subspace of $\mathcal{H}$.
To $\mathcal{K}$ we can extend the new norm $\norm{\cdot}_\sharp$ as follows.
For $\zeta \in \mathcal{K}$ define
\begin{equation*}
  \norm{\zeta}_\sharp
  = \mathrm{inf}
  \defset{
    \mathrm{max}\left\{ \norm{\zeta_1}_\sharp, \norm{\zeta_2}_\sharp \right\}
    }{\zeta_1,\zeta_2 \in \mathcal{P}^\sharp_b, \zeta_1-\zeta_2 = \zeta}.
\end{equation*}

It is easily seen that if $z \in \mathcal{M}_{sa}$ corresponds to
$\zeta \in \mathcal{K}$, we have
\begin{equation*}
  \mathrm{max}\left\{\norm{z_+}, \norm{z_-}\right\} = \norm{z}
   = \norm{\zeta}_\sharp
   = \mathrm{max}\left\{\norm{\zeta_+}_\sharp, \norm{\zeta_-}_\sharp \right\}.
\end{equation*}

\section{Jordan structure on $\mathcal{K}+i\mathcal{K}$}
First we define the square operation for vectors in $\mathcal{K}$.

\begin{definition}
If $\zeta$ is a real linear combination of mutually operationally orthogonal
projective vectors, i.e. $\zeta = \sum_k c_k \zeta_k$ where $c_k \in \mathbb{R}$
and $\{\zeta_k\}$ are mutually operationally orthogonal, then we define the
square of $\zeta$ as follows.
\begin{equation*}
  \zeta^2 = \sum_k c_k^2 \zeta_k.
\end{equation*}
\end{definition}

As we have seen in Lemma \ref{fund}, mutually operationally
orthogonal projective vectors $\{\zeta_k\}$ correspond to mutually orthogonal
projections $\{e_k\}$. Thus the square of a real linear combination
$\sum_k c_k e_k$ equals $\sum_k c_k^2 e_k$ and for these vectors the definition
of square is consistent.

The set of vectors which are real linear combinations of mutually operationally
orthogonal projective vectors is dense in $\mathcal{K}$ in the sharp norm
defined in Section \ref{cone}.
In fact, these vectors correspond to real linear combinations of mutually orthogonal projections in $\mathcal{M}$, i.e. self-adjoint operators with
finite spectra. 

Since the sharp norm on $\mathcal{K}$ is consistent with the
operator norm on $\mathcal{M}$, we can extend the definition of square
to $\mathcal{K}$ by continuity. We have the following.
\begin{equation*}
  \mbox{If } \zeta = z\xi_0 \mbox{ for } z \in \mathcal{M}_{sa} \mbox{, then }
  \zeta^2 = z^2\xi_0.
\end{equation*}

Once we have defined the square operation on $\mathcal{K}$, we can define
Jordan polynomials as follows.
For $\eta$ and $\zeta$ in $\mathcal{K}$ let us define
\begin{equation*}
  \eta\zeta + \zeta\eta = (\eta+\zeta)^2 - \eta^2 - \zeta^2.
\end{equation*}
Using this, for $\zeta = \zeta_1 + i\zeta_2 \in \mathcal{K}+i\mathcal{K}$
we put
\begin{equation*}
  \zeta^2 = \zeta_1^2 + i(\zeta_1\zeta_2+\zeta_2\zeta_2) - \zeta_2^2.
\end{equation*}
As for vectors in $\mathcal{K}$, we define the ``Jordan product'' on
$\mathcal{K} +i\mathcal{K}$ by
\begin{equation*}
  \eta\zeta + \zeta\eta = (\eta+\zeta)^2 - \eta^2 - \zeta^2.
\end{equation*}
Using this, finally we define
\begin{equation*}
  \zeta\eta\zeta =
  \frac{1}{2}\left[(\zeta\eta+\eta\zeta)\zeta+\zeta(\zeta\eta+\eta\zeta)\right]
   -\frac{1}{2}\left(\zeta^2\eta+\eta\zeta^2\right).
\end{equation*}

If $\eta = y\xi_0$ and $\zeta = z\xi_0$ for $y,z \in \mathcal{M}$, then
it follows that $\zeta\eta\zeta = zyz\xi_0$. This follows because we have
defined square and Jordan polynomials on $\mathcal{K}$ consistently.

If we fix $\zeta$, we give names to the following mappings.
\begin{align*}
  \mathrm{c}_{\zeta}&:\mathcal{K}+i\mathcal{K} \ni \eta \longmapsto
    \zeta\eta\zeta \in \mathcal{K}+i\mathcal{K}, \\
  \mathrm{od}_{\zeta}&:\mathcal{K}+i\mathcal{K} \ni \eta \longmapsto
    \eta - \corner{\zeta}{\eta} - \corner{\zeta^\perp}{\eta}
    \in \mathcal{K}+i\mathcal{K}.
\end{align*}

 Let $\eta = y\xi_0$ and $\zeta = e\xi_0$ where $e$ is a projection.
Then we see that
\begin{align*}
\corner{\zeta}{\eta} &= eye\xi_0, \mbox{ and} \\
\offdiag{\zeta}{\eta} &= y\xi_0 - eye\xi_0 - e^\perp ye^\perp\xi_0
 = \left[eye^\perp + e^\perp ye\right]\xi_0
\end{align*}
correspond to the corner of $y$ and the off-diagonal part of $y$,
respectively.

\section{Recovery of projections in $\mathcal{M}$ in the case when
$\mathcal{M}^\sigma = \mathbb{C}I$}

Let $p$ be a projection of $\mathcal{B}(\mathcal{H})$. We seek a
necessary and sufficient condition for $p$ to be a projection of $\mathcal{M}$.

We need a criterion for a projection in $\mathcal{M}$ to be fixed by
the modular automorphism.
\begin{lemma}\label{fixed}
Let $e$ be a projection
in $\mathcal{M}$. If $px\xi_0 = xe\xi_0$ holds for all $x \in \mathcal{M}$,
then we have $e \in \mathcal{M}^\sigma$ and $p = JeJ$.
\begin{proof}
Note that we get $p\xi_0 = e\xi_0$ if we use the assumption with $x = I$.

Again by the assumption it follows that
\begin{eqnarray*}
 \langle xe\xi_0, \xi_0\rangle &=& \langle px\xi_0, \xi_0\rangle \\
                               &=& \langle x\xi_0, p\xi_0\rangle \\
                               &=& \langle x\xi_0, e\xi_0\rangle \\
                               &=& \langle ex\xi_0, \xi_0\rangle.
\end{eqnarray*}
This implies that $e \in \mathcal{M}^\sigma$\cite{sz}. In particular, we have
\begin{equation*}
  e\xi_0 = Se\xi_0 = J\Delta^\frac{1}{2}e\xi_0 = Je\xi_0.
\end{equation*}

Now the equality $JeJx\xi_0 = xJeJ\xi_0 = xe\xi_0 = px\xi_0$
and the cyclicity of $\xi_0$ complete the proof. 
\end{proof}
\end{lemma}

Recall that $S = J\Delta^\frac{1}{2}$ can be defined in terms of
$\overline{\mathcal{K}}$ \cite{rvd}.

\begin{theorem}\label{recover}
Let $p$ be a projection in $\mathcal{B}(\mathcal{H})$.
There is a projection $e \in \mathcal{M}$ and a central projection
$q \in \mathcal{M}$ such that $q^\perp e \in \mathcal{M}^\sigma$ and
$p = qe + Jq^\perp eJ$ if and only if
the following hold:
\begin{enumerate}
\item $p\xi_0 \le \xi_0$.
\item If $\zeta \le p\xi_0$, then $p\zeta = \zeta$. \label{cnp1}
\item If $\zeta \le p^\perp\xi_0$, then $p^\perp\zeta = \zeta$. \label{cnp2}
\item For every vector $\xi \in \mathcal{K}+i\mathcal{K}$ we have
 $p\xi \in \mathcal{K}+i\mathcal{K}$ and
 \begin{enumerate}
 \item $\corner{p\xi_0}{p\ \offdiag{p\xi_0}{\xi}} = 0$, \label{cn1}
 \item $\corner{p^\perp\xi_0}{p\ \offdiag{p\xi_0}{\xi}} = 0$, \label{cn2}
 \item $\left(p\ \offdiag{p\xi_0}{\xi}\right)^2 = 0$, \label{sq1}
 \item $\left(p^\perp\ \offdiag{p\xi_0}{\xi}\right)^2 = 0$, \label{sq2}
 \item $Sp\ \offdiag{p\xi_0}{\xi} = p^\perp S\ \offdiag{p\xi_0}{\xi}$.
    \label{star}
 \end{enumerate} \label{remain}
\end{enumerate}

\begin{proof}
First let us show the ``only if'' part.
In this case, we have
\begin{equation*}
  p\xi_0 = qe\xi_0 + Jq^\perp eJ\xi_0 = qe\xi_0 + q^\perp e\xi_0 = e\xi_0
   \le \xi_0,
\end{equation*}
hence the first part of the conditions is satisfied. For the second condition,
if $\zeta = z\xi_0 \le p\xi_0 = e\xi_0$, then the support of $z$ is less than or
equal to $e$ and we have
\begin{equation*}
  p\zeta = qez\xi_0 + zJeq^\perp J\xi_0 = qez\xi_0 + zeq^\perp \xi_0 = z\xi_0
   = \zeta.
\end{equation*}
Similar proof works for the third. To see the conditions of the fourth,
let $\xi = x\xi_0 \in \mathcal{K}+i\mathcal{K}$. We note that
\begin{align*}
\corner{p\xi_0}{\xi} &= \corner{e\xi_0}{x\xi_0} = exe\xi_0, \\
\offdiag{p\xi_0}{\xi} &= \offdiag{e\xi_0}{x\xi_0}
 = \left[exe^\perp + e^\perp xe \right]\xi_0, \\
p\ \offdiag{p\xi_0}{\xi}
 &= \left[qexe^\perp + q^\perp e^\perp xe\right]\xi_0, \\
p^\perp \ \offdiag{p\xi_0}{\xi}
 &= \left[qe^\perp xe + q^\perp exe^\perp \right]\xi_0, \\
Sp\ \offdiag{p\xi_0}{\xi}
 &= \left[qe^\perp x^*e + q^\perp ex^*e^\perp \right]\xi_0, \\
p^\perp S\ \offdiag{p\xi_0}{\xi}
 &= (qe^\perp+Jq^\perp e^\perp J)\left[e^\perp x^*e+ex^*e^\perp \right]\xi_0 \\
 &= \left[qe^\perp x^*e + q^\perp ex^*e^\perp \right]\xi_0. 
\end{align*}
Thus it is easy to see that each of the conditions is valid.

We turn to the ``if'' part. Let $p$ satisfy the conditions of the statement. 

Take $x \in \mathcal{M}$ satisfying
$x = exe^\perp$. If we use the matrix, $x$ takes the following form.
\begin{equation*}
\bordermatrix{
                & \mathrm{Ran}(e) & \mathrm{Ran} (e^\perp) \cr
 \mathrm{Ran}(e)&              0 &                    X \cr
 \mathrm{Ran}(e^\perp)&        0 &                    0
}.
\end{equation*}
Then it holds that $\offdiag{p\xi_0}{x\xi_0} = x\xi_0$. 

By assumption \ref{remain}, there exists $y \in \mathcal{M}$ such that
$px\xi_0 = y\xi_0$. In addition, by assumptions
\ref{cn1} and \ref{cn2}, we have $eye = e^\perp ye^\perp = 0$,
i.e. $y$ has trivial corners.
By assumption \ref{sq1}, it follows $y^2 = 0.$ Hence $y$ takes
the following form.
\begin{equation*}
  y = \left(
       \begin{array}{cc}
         \scalebox{1.5}{0} & \begin{array}{ccc}
               y_1 & 0 & 0 \\
                 0 & 0 & 0 \\
                 0 & 0 & 0
             \end{array} \\
         \begin{array}{ccc}
             0 &   0 & 0 \\
             0 & y_2 & 0 \\
             0 &   0 & 0
         \end{array} & \scalebox{1.5}{0}
       \end{array} \right),
\end{equation*}
where we decomposed $\mathrm{Ran}(e)$ and $\mathrm{Ran}(e^\perp)$ as follows.
\begin{align*}
\mathrm{Ran}(e) &= \mathrm{Dom}(e^\perp ye) \oplus \mathrm{Ran}(eye^\perp)
 \oplus \left(\mathrm{Ran}(e) \ominus
               \mathrm{Dom}(e^\perp ye) \ominus
               \mathrm{Ran}(eye^\perp)\right), \\
\mathrm{Ran}(e^\perp) &= \mathrm{Dom}(eye^\perp) \oplus \mathrm{Ran}(e^\perp ye)
  \oplus \left(\mathrm{Ran}(e^\perp) \ominus
               \mathrm{Dom}(eye^\perp) \ominus
               \mathrm{Ran}(e^\perp ye) \right).
\end{align*}
Subspaces which appear here are mutually orthogonal because
the square of $y$ vanishes.

According to this, we further decompose $x$.
\begin{equation*}
  x = \left(
       \begin{array}{cc}
         \scalebox{1.5}{0} & \begin{array}{ccc}
               x_1 & x_2 & x_3 \\
               x_4 & x_5 & x_6 \\
               x_7 & x_8 & x_9 
             \end{array} \\
         \begin{array}{ccc}
             0 &   0 & 0 \\
             0 &   0 & 0 \\
             0 &   0 & 0
         \end{array} & \scalebox{1.5}{0}
       \end{array} \right).
\end{equation*}
By assumption \ref{sq2}, the square of $p^\perp x\xi_0 = (x-y)\xi_0$
must vanish.
\begin{align*}
  x-y &= \left(
       \begin{array}{cc}
         \scalebox{1.5}{0} & \begin{array}{ccc}
               x_1-y_1 & x_2 & x_3 \\
               x_4     & x_5 & x_6 \\
               x_7     & x_8 & x_9 
             \end{array} \\
         \begin{array}{ccc}
             0 &      0 & 0 \\
             0 &   -y_2 & 0 \\
             0 &      0 & 0
         \end{array} & \scalebox{1.5}{0}
       \end{array} \right), \\
 (x-y)^2 &= \left(
       \begin{array}{cc}
         \begin{array}{ccc}
               0 & -x_2 y_2 & 0 \\
               0 & -x_5 y_2 & 0 \\
               0 & -x_8 y_2 & 0 
         \end{array} & \scalebox{1.5}{0} \\
         \scalebox{1.5}{0} & \begin{array}{ccc}
                                   0 &          0 &        0 \\
                            -y_2 x_4 &   -y_2 x_5 & -y_2 x_6 \\
                                   0 &            0 &      0
         \end{array}
       \end{array} \right).
\end{align*}
Then it follows that $x_2 = x_4 = x_5 = x_6 = x_8 = 0$.

If we use assumption \ref{star}, then we get
\begin{equation*}
  px^*\xi_0 = pSx\xi_0 = Sp^\perp x\xi_0 = (x^*-y^*)\xi_0.
\end{equation*}
Applying assumption \ref{sq1} to $\xi = (x+x^*)\xi_0$, the square of $p(x+x^*)\xi_0 = (y+x^*-y^*)\xi_0$ vanishes. 
\begin{align*}
 & y+x^*-y^* = \left(
       \begin{array}{cc}
         \scalebox{1.5}{0} & \begin{array}{ccc}
                               y_1 &      0 & 0 \\
                                 0 & -y_2^* & 0 \\
                                 0 &      0 & 0
             \end{array} \\
         \begin{array}{ccc}
             x_1^*-y_1^* &   0 & x_7^* \\
                       0 & y_2 &     0 \\
                   x_3^* &   0 & x_9^*
         \end{array} & \scalebox{1.5}{0}
       \end{array} \right), \\
& (y+x^*-y^*)^2 \\
&   =  \begin{pmatrix}
         \begin{array}{ccc}
             y_1(x_1^*-y_1^*) &         0 & y_1x_7^* \\
                            0 & -y_2^*y_2 &        0 \\
                            0 &         0 &        0 
         \end{array} & \scalebox{1.5}{0} \\
         \scalebox{1.5}{0} & \begin{array}{ccc}
                                (x_1^*-y_1^*)y_1 &         0 & 0 \\
                                               0 & -y_2y_2^* & 0 \\
                                        x_3^*y_1 &         0 & 0
                             \end{array}
       \end{pmatrix}.
\end{align*}

Thus it follows that $y_2 = x_3 = x_7 = 0$ and $x_1 = y_1$.

Summing up, for every $x = exe^\perp \in \mathcal{M}$ we have
\begin{align*}
  x &= \left(
       \begin{array}{cc}
         \scalebox{1.5}{0} & \begin{array}{ccc}
                               x_1 & 0 &   0 \\
                                 0 & 0 &   0 \\
                                 0 & 0 & x_9
                             \end{array} \\
         \begin{array}{ccc}
             0 & 0 & 0 \\
             0 & 0 & 0 \\
             0 & 0 & 0
         \end{array} & \scalebox{1.5}{0}
       \end{array} \right), \\
  y\xi_0 = px\xi_0 &= \left(
       \begin{array}{cc}
         \scalebox{1.5}{0} & \begin{array}{ccc}
                               x_1 & 0 & 0 \\
                                 0 & 0 & 0 \\
                                 0 & 0 & 0
                             \end{array} \\
         \begin{array}{ccc}
             0 & 0 & 0 \\
             0 & 0 & 0 \\
             0 & 0 & 0
         \end{array} & \scalebox{1.5}{0}
       \end{array} \right)\xi_0.
\end{align*}
The point is that $\mathrm{Dom}(y)$ and $\mathrm{Dom}(x-y)$,
$\mathrm{Ran}(y)$ and $\mathrm{Ran}(x-y)$ are mutually orthogonal,
respectively.

If we take another element $z = eze^\perp \in \mathcal{M}$ and put
$w\xi_0 = pz\xi_0$, then by the same argument we see that
$\mathrm{Dom}(w)$ and $\mathrm{Dom}(z-w)$,
$\mathrm{Ran}(w)$ and $\mathrm{Ran}(z-w)$ are mutually orthogonal,
respectively.
In addition, by noting that $w+x-y = e(w+x-y)e^\perp$ and
$p(w+x-y)\xi_0 = w\xi_0$,
it follows that $\mathrm{Dom}(x-y) \perp \mathrm{Dom}(w)$
and $\mathrm{Ran}(x-y) \perp \mathrm{Ran}(w)$. Similarly it holds
that $\mathrm{Dom}(z-w) \perp \mathrm{Dom}(y)$
and $\mathrm{Ran}(z-w) \perp \mathrm{Ran}(y)$.
Then let us define $f_1$ (resp. $f_3$) to be the projection onto the supremum of
such $\mathrm{Ran}(x-y)$'s (resp. $\mathrm{Dom}(x-y)$'s) where
$x = exe^\perp$ runs all the elements of this form in $\mathcal{M}$
and put $f_2 = e - f_1$, $f_4 = e^\perp - f_3$. They are mutually orthogonal
projections of
$\mathcal{M}$.

Using them every $x = exe^\perp \in \mathcal{M}$ is decomposed
as follows.
\begin{equation*}
 \bordermatrix{
   & \mathrm{Ran}(f_1) & \mathrm{Ran}(f_2) & \mathrm{Ran}(f_3) & \mathrm{Ran}(f_4) \cr
   \mathrm{Ran}(f_1) & 0 & 0 & x_1 &   0 \cr
   \mathrm{Ran}(f_2) & 0 & 0 &   0 & x_2 \cr
   \mathrm{Ran}(f_3) & 0 & 0 &   0 &   0 \cr
   \mathrm{Ran}(f_4) & 0 & 0 &   0 &   0 
 }.
\end{equation*}
According to this decomposition, it is easy to see that
every $x \in \mathcal{M}$ must have the following form.
\begin{equation*}
  x = \left( \begin{array}{cccc}
        x_1 & 0 & x_3 &   0 \\
        0 & x_2 &   0 & x_4 \\
        x_5 & 0 & x_7 &   0 \\
        0 & x_6 &   0 & x_8
      \end{array}\right).
\end{equation*}
Put $q = f_1 + f_3$. This is clearly a central projection.

Since $p$ preserves vectors of the set
$\defset{\zeta}{\zeta \le p\xi_0 = e\xi_0}$
by assumption \ref{cnp1},
it holds that $p\ exe\xi_0 = exe\xi_0$ for $x \in \mathcal{M}$.
Similarly, by assumption \ref{cnp2}, we see
$p^\perp \ e^\perp x e^\perp \xi_0 = e^\perp x e^\perp \xi_0$, 
hence $p\ e^\perp x e^\perp \xi_0 = 0$.

Now, letting $x$ be an arbitrary element of $\mathcal{M}$, $p$
acts on $x\xi_0$ as follows.
\begin{align*}
px\xi_0 &= p\left( \begin{array}{cccc}
        x_1 & 0 & x_3 &   0 \\
        0 & x_2 &   0 & x_4 \\
        x_5 & 0 & x_7 &   0 \\
        0 & x_6 &   0 & x_8
      \end{array}\right)\xi_0
      = \left( \begin{array}{cccc}
        x_1 & 0 & x_3 & 0 \\
        0 & x_2 &   0 & 0 \\
        0 &   0 &   0 & 0 \\
        0 & x_6 &   0 & 0
      \end{array}\right)\xi_0 \\
      &= (qex + q^\perp xe)\xi_0.
\end{align*}
Then using the cyclicity of $\xi_0$ and Lemma \ref{fixed},
we arrive at the conclusion that $p = qe + Jq^\perp eJ$. 
\end{proof}
\end{theorem}

\begin{corollary}
If $\mathcal{M}^\sigma = \mathbb{C}I$, then the conditions in
Theorem \ref{recover} assure that $p$ is a projection of
$\mathcal{M}$.
\end{corollary}

\subsubsection*{Acknowledgements.}
I am truly grateful to my supervisor Yasuyuki Kawahigashi for
his helpful comments and supports. I also would like to
thank Roberto Longo and Yasuhide Miura for their valuable advice.


\begin{thebibliography}{10}

\bibitem{araki}
H.~Araki.
\newblock Some properties of modular conjugation operator of von {N}eumann
  algebras and a non-commutative {R}adon-{N}ikodym theorem with a chain rule.
\newblock {\em Pacific J. Math.}, 50:309--354, 1974.

\bibitem{az}
H.~Araki and L.~Zsid{\'o}.
\newblock Extension of the structure theorem of {B}orchers and its application
  to half-sided modular inclusions.
\newblock {\em Rev. Math. Phys.}, 17(5):491--543, 2005.

\bibitem{baumgaertel}
H.~Baumg{\"a}rtel.
\newblock {\em Operatoralgebraic methods in quantum field theory}.
\newblock Akademie Verlag, Berlin, 1995.

\bibitem{borchers}
H.-J. Borchers.
\newblock Half-sided modular inclusions and structure analysis in quantum field
  theory.
\newblock In {\em Operator algebras and quantum field theory (Rome, 1996)},
  pages 589--608. Int. Press, Cambridge, MA, 1997.

\bibitem{connes1}
A.~Connes.
\newblock Groupe modulaire d'une alg\`ebre de von {N}eumann.
\newblock {\em C. R. Acad. Sci. Paris S\'er. A-B}, 274:A1923--A1926, 1972.

\bibitem{connes2}
A.~Connes.
\newblock Caract\'erisation des espaces vectoriels ordonn\'es sous-jacents aux
  alg\`ebres de von {N}eumann.
\newblock {\em Ann. Inst. Fourier (Grenoble)}, 24(4):x, 121--155 (1975), 1974.

\bibitem{davidson}
D.~R. Davidson.
\newblock Endomorphism semigroups and lightlike translations.
\newblock {\em Lett. Math. Phys.}, 38(1):77--90, 1996.

\bibitem{jr}
N.~Jacobson and C.~E. Rickart.
\newblock Jordan homomorphisms of rings.
\newblock {\em Trans. Amer. Math. Soc.}, 69:479--502, 1950.

\bibitem{kadison}
R.~V. Kadison.
\newblock Isometries of operator algebras.
\newblock {\em Ann. of Math. (2)}, 54:325--338, 1951.

\bibitem{rvd}
M.~A. Rieffel and A.~van Daele.
\newblock A bounded operator approach to {T}omita-{T}akesaki theory.
\newblock {\em Pacific J. Math.}, 69(1):187--221, 1977.

\bibitem{sz}
{\c{S}}.~Str{\u{a}}til{\u{a}} and L.~Zsid{\'o}.
\newblock {\em Lectures on von {N}eumann algebras}.
\newblock Editura Academiei, Bucharest, 1979.

\bibitem{wiesbrock}
H.-W. Wiesbrock.
\newblock Half-sided modular inclusions of von-{N}eumann-algebras.
\newblock {\em Comm. Math. Phys.}, 157(1):83--92, 1993.

\end{thebibliography}
\end{document}